\documentclass[12pt]{amsart}

\usepackage{amsmath,amssymb,amsthm,enumerate}

\begin{document}
\numberwithin{equation}{section}

\def\1#1{\overline{#1}}
\def\2#1{\widetilde{#1}}
\def\3#1{\widehat{#1}}
\def\4#1{\mathbb{#1}}
\def\5#1{\frak{#1}}
\def\6#1{{\mathcal{#1}}}

\def\C{{\4C}}
\def\R{{\4R}}
\def\N{{\4N}}
\def\Z{{\4Z}}

\def \im{\text{\rm Im }}
\def \re{\text{\rm Re }}
\def \Char{\text{\rm Char }}
\def \supp{\text{\rm supp }}
\def \codim{\text{\rm codim }}
\def \Ht{\text{\rm ht }}
\def \Dt{\text{\rm dt }}
\def \hO{\widehat{\mathcal O}}
\def \cl{\text{\rm cl }}
\def \bR{\mathbb R}
\def \bC{\mathbb C}
\def \bP{\mathbb P}
\def \C{\mathbb C}
\def \bL{\mathbb L}
\def \bZ{\mathbb Z}
\def \bN{\mathbb N}
\def \scrF{\mathcal F}
\def \scrK{\mathcal K}
\def \scrM{\mathcal M}
\def \cR{\mathcal R}
\def \scrJ{\mathcal J}
\def \scrA{\mathcal A}
\def \scrO{\mathcal O}
\def \scrV{\mathcal V}
\def \scrL{\mathcal L}
\def \scrE{\mathcal E}

\title[Normal form for 1-infinite type hypersurfaces]{A normal form for 1-infinite type hypersurfaces in $\bC^2$. I. Formal Theory.}
\author{Peter Ebenfelt}
\address{Department of Mathematics, University of California at San Diego, La Jolla, CA 92093-0112}
\email{pebenfel@math.ucsd.edu}

\author{Bernhard Lamel}
\address{Department of Mathematics, University of Vienna, Vienna, Austria}
\email{bernhard.lamel@univie.ac.at}

\author[D. Zaitsev]{Dmitri Zaitsev}
\address{D. Zaitsev: School of Mathematics, Trinity College Dublin, Dublin 2, Ireland}
\email{zaitsev@maths.tcd.ie}

\thanks{The first author was supported in part by the NSF grant DMS-1301282. The second author
was supported by the Austrian Federal Ministry of Research through START prize Y377.
The third author was supported by the Science Foundation Ireland grant 10/RFP/MTH2878.}

\thanks{2010 {\em Mathematics Subject Classification}. 32H02, 32V40}

\begin{abstract} In this paper, we study real hypersurfaces $M$ in $\bC^2$ at points $p\in M$ of infinite type. The degeneracy of $M$ at $p$ is assumed to be the least possible, namely such that the Levi form vanishes to first order in the CR transversal direction. A new phenomenon, compared to known normal forms in other cases, is the presence of resonances as roots of a universal polynomial in the $7$-jet of the defining function of $M$. The main result is a complete (formal) normal form at points $p$ with no resonances.
Remarkably, our normal form at such infinite type points resembles closely the Chern-Moser normal form at Levi-nondegenerate points.
For a fixed hypersurface, its normal forms are parametrized by $S^1\times \bR^*$, and as a corollary we find that the automorphisms in the stability group of $M$ at $p$ without resonances are determined by their $1$-jets at $p$. In the last section, as a contrast, we also give examples of hypersurfaces with arbitrarily high resonances that possess families of distinct automorphisms whose jets agree up to the resonant order.
\end{abstract}

\maketitle   

\def\Label#1{\label{#1}}


\def\cn{{\C^n}}
\def\cnn{{\C^{n'}}}
\def\ocn{\2{\C^n}}
\def\ocnn{\2{\C^{n'}}}


\def\dist{{\rm dist}}
\def\const{{\rm const}}
\def\rk{{\rm rank\,}}
\def\id{{\sf id}}
\def\aut{{\sf aut}}
\def\Aut{{\sf Aut}}
\def\CR{{\rm CR}}
\def\GL{{\sf GL}}
\def\Re{{\sf Re}\,}
\def\Im{{\sf Im}\,}
\def\span{\text{\rm span}}

\def\codim{{\rm codim}}
\def\crd{\dim_{{\rm CR}}}
\def\crc{{\rm codim_{CR}}}

\def\phi{\varphi}
\def\eps{\varepsilon}
\def\d{\partial}
\def\a{\alpha}
\def\b{\beta}
\def\g{\gamma}
\def\G{\Gamma}
\def\D{\Delta}
\def\Om{\Omega}
\def\k{\kappa}
\def\l{\lambda}
\def\L{\Lambda}
\def\z{{\bar z}}
\def\w{{\bar w}}
\def\Z{{\1Z}}
\def\t{\tau}
\def\th{\theta}

\emergencystretch15pt
\frenchspacing

\newtheorem{Thm}{Theorem}[section]
\newtheorem{Cor}[Thm]{Corollary}
\newtheorem{Pro}[Thm]{Proposition}
\newtheorem{Lem}[Thm]{Lemma}

\theoremstyle{definition}\newtheorem{Def}[Thm]{Definition}

\theoremstyle{remark}
\newtheorem{Rem}[Thm]{Remark}
\newtheorem{Exa}[Thm]{Example}
\newtheorem{Exs}[Thm]{Examples}

\def\bl{\begin{Lem}}
\def\el{\end{Lem}}
\def\bp{\begin{Pro}}
\def\ep{\end{Pro}}
\def\bt{\begin{Thm}}
\def\et{\end{Thm}}
\def\bc{\begin{Cor}}
\def\ec{\end{Cor}}
\def\bd{\begin{Def}}
\def\ed{\end{Def}}
\def\br{\begin{Rem}}
\def\er{\end{Rem}}
\def\be{\begin{Exa}}
\def\ee{\end{Exa}}
\def\bpf{\begin{proof}}
\def\epf{\end{proof}}
\def\ben{\begin{enumerate}}
\def\een{\end{enumerate}}

\section{Introduction}

The main objective in this paper is to construct a normal form,
unique up to the action of a finite dimensional group,
for a class of real hypersurfaces $M$ in $\bC^2$ at a point $p$ of {\em infinite type}; recall that $p$ is of infinite type if (and only if, in the real-analytic case) there is a complex curve through $p$ contained in $M$. The new phenomenon here, compared to previously known cases of normal forms for CR manifolds,
is the presence of so-called {\em resonances}.
A resonance is an integral root of a polynomial, called the {\em characteristic polynomial}, whose coefficients
are polynomials in the $7$-jet of the defining equation of the hypersurface.
If a hypersurface has no resonances, we obtain a normal form
unique up to rotations and scaling.
If, on the other hand, resonances are present,
the same normalization conditions are obtained for
all terms except the resonant ones (of which there are always at most finitely many).

There is a vast literature on normal forms of real hypersurfaces in $\bC^{n+1}$ at points of finite (commutator) type, starting with the seminal paper \cite{CM74} from 1974 by S.-S. Chern and J. Moser providing a normal form for Levi nondegenerate hypersurfaces. The Chern-Moser normal form is {\it convergent} in the sense that if $M$ is a real-analytic hypersurface, then the transformation to normal form is holomorphic (given by a convergent formal transformation) and the resulting equation in normal form converges to a defining equation for the transformed hypersurface. Following the lead of Chern-Moser, numerous authors have constructed normal forms for various classes of real hypersurfaces at points of finite type; we mention here only \cite{Wong82}, \cite{Stanton91}, \cite{E98a}, \cite{E98b}, \cite{EHZ05}, \cite{Kolar05}, \cite{KMZ14}, \cite{KoZ14a}, \cite{KoZ14b}. These normal forms are all formal (not known to be convergent, or in some cases even known to not be convergent \cite{Kolar12}), with the exception of the very recent \cite{KoZ14a}, \cite{KoZ14b}. The normal form we construct in this paper is formal. We should point out, however, that there are general results concerning convergence of formal invertible mappings between real-analytic CR manifolds (see \cite{BER00}, \cite{BMR02}) that apply in the finite type situations treated in the papers mentioned above. As a consequence, questions about biholomorphic mappings (such as, e.g., their existence) between real-analytic CR manifolds of finite type (that are also holomorphically nondegenerate; \cite{BMR02}) can often be reduced to the analogous questions about formal mappings, and for the latter it suffices that the manifolds are in formal normal form. For the class of infinite type hypersurfaces considered in this paper, the corresponding convergence result for formal mappings between real-analytic hypersurfaces is known as well (\cite{JuhlinLamel13}; cf.\ also the unpublished thesis \cite{Juhlin07}).

As mentioned, there is a vast literature on normal forms for real hypersurfaces at points of finite type, but the normal form presented here is (to the best of the authors' knowledge) the first systematic result of this kind at points of infinite type. There is, however, a previous paper by the authors \cite{ELZ09}, in which new invariants are introduced for real hypersurfaces in $\bC^2$ and a (formal) normal form is constructed for a certain class of hypersurfaces identified by conditions on these invariants. The class so identified contains some hypersurfaces of infinite type, but is in fact completely disjoint from the class considered in this paper. The main objective in \cite{ELZ09} was to provide conditions in terms of the new invariants that would guarantee triviality (or discreteness) of the stability group of the hypersurface.
The normal form in that paper was {\it ad hoc} and its main purpose was a means to prove the result about the stability group. There are also the results in \cite{KL14}, in which a dimension bound was proved by means of an ``abstract'' normal form construction (which however does not produce a normal form at all in the usual sense).

The main objective in this paper is the construction of a normal form, which is modeled on the Chern-Moser normal form in the Levi-nondegenerate case, for a rich and natural class of infinite type hypersurfaces. For a comparison of the normal form in this paper (described below) and the Chern-Moser normal form, we recall, for the reader's convenience, the Chern-Moser normal form for a smooth Levi nondegenerate hypersurface $M$ through $0$ in $\bC^2$: There are formal holomorphic  coordinates $(z,w)$ near $0\in \bC^2$ such that $M$ is given locally by
\begin{equation}
\Im w=\phi(z,\bar z,\Re w),
\end{equation}
where the (Hermitian) formal power series $\phi(z,\bar z,u)$ is of the form
\begin{equation}
\phi(z,\bar z,u)=|z|^2+\sum_{a,b\geq 0}N_{ab}(u)z^a\bar z^b,\quad \left(N_{ab}(u)=\overline{N_{ba}(u)},\ u\in \bR\right),
\end{equation}
satisfying the following normalization conditions
\begin{equation}
N_{a0}(u)=N_{a1}(u)=0\quad a\geq 0,
\end{equation}
and
\begin{equation}
N_{22}(u)=N_{32}(u)=N_{33}(u)=0.
\end{equation}
This normal form is unique modulo the action of the stability group of the sphere ($(\Im w=|z|^2$ in these coordinates).
Our normal form in Theorem \ref{main0} below is very similar.

\subsection{The main result}
We shall now describe our main result more precisely. Let $M\subset \bC^2$ be a smooth real hypersurface with $p\in M$. After an affine linear tranformation, we find local holomorphic coordinates $(z,w)$, vanishing at $p$, such that the real tangent space to $M$ at $0$ is spanned by
 $\Re \partial/\partial z,\Im \partial/\partial z, \Re \partial/\partial w$, and $M$ is given locally as a graph
\begin{equation}\Label{defeq0}
\Im w=\phi(z,\bar z, \Re w),
\end{equation}
where $\phi(0)=0$ and $d\phi(0)=0$. We shall assume the following:
\begin{itemize}
\item [(1)] $M$ is of {\it infinite type} at $p=0$, i.e., for any $m$, there is a holomorphic curve $\Gamma_m\colon \bC\to \bC^2$ with $\Gamma_m(0)=0$, which is tangent to $M$ to order $m$ at $0$.
\item[(2)] There is a smooth curve $\gamma\colon (-\epsilon,\epsilon)\to M$ with $\gamma(0)=0$ (necessarily transverse to the complex tangent space of $M$ at $0$) such that the Levi form of $M$ vanishes to first order at $0$ along $\gamma$, i.e., $(\mathcal L\circ \gamma)'(0)\neq 0$, where $\mathcal L$ is any representative of the Levi form of $M$.
\end{itemize}
In view of (1) and (2), we can assume
\begin{equation}\Label{order}
\phi(z,\bar z, u) = z\bar z u + O(|(z,\bar z,u)|^4).
\end{equation}
We shall then introduce a monic polynomial $P(k, j_0^7\phi)$ in $k\in\C$ and the 7-jet of $\phi$ at $0$, see Definition \ref{generalpos},
and the following {\em nonresonant condition}:
\begin{itemize}
\item[(3)] $P(k, j_0^7\phi)$ has no integral roots $k\ge 2$ (which we call {\em resonances}).
\end{itemize}
We mention here that (3) holds for $j^7_0\phi$ in a specific open and dense subset $\Omega$ of $J^7_0(\bC\times\bR)$, the space of 7-jets at 0 of smooth functions in $\bC\times\bR$.
Indeed, since $P$ is monic in $k$, the set $\Omega$ is locally
determined by finitely many polynomial inequalities
(for a finite set of possible roots $k$).

Our main result is a formal normal form for the hypersurface $M$ at $p=0$. This normal form is unique, as is the formal transformation to normal form, modulo the action of the 2-dimensional group $S^1\times \bR^*$, where $S^1$ denotes the unit circle and $\bR^*:=\bR\setminus \{0\}$. Since our normal form is formal, we shall formulate our result for formal hypersurfaces.

For our purposes, a {\it formal hypersurface} through $0$ in the coordinates $(z,w)\in \bC^2$ is an object associated to a graph equation of the form \eqref{defeq0}, where $\phi(z,\bar z,u)$ is a formal power series in $z,\bar z,u$ such that $\phi(0)=0$. Clearly, a smooth hypersurface $M$ through $p=0$ as above, defines a formal hypersurface via the Taylor series of the smooth graphing function $\phi$ in \eqref{defeq0}; by an abuse of notation, we shall continue to denote the formal hypersurface by $M$ and the formal graphing power series by $\phi(z,\bar z,u)$. We note that two distinct smooth hypersurfaces through $p=0$ may define the same formal hypersurface; this happens if and only if the two hypersurfaces are tangent to infinite order at $0$. We also note that a smooth hypersurface $M$ satisfies Conditions (1) and (2) above if and only if its associated formal hypersurface $M$ satisfies the following conditions:
\begin{itemize}
\item [(1')] $M$ is of {\it infinite type} at $p=0$; i.e., there is a formal holomorphic curve $\Gamma\colon \bC\to \bC^2$ with $\Gamma(0)=0$, which is contained in $M$ (formally); see \cite{BER99a}.
\item[(2')] There is a formal curve $\gamma\colon \bR\to M$ with $\gamma(0)=0$ (necessarily transverse to the complex tangent space of $M$ at $0$) such that the Levi form $\mathcal L$ of $M$ satisfies $(\mathcal L\circ \gamma)'(0)\neq 0$.
\end{itemize}
Also, note that Condition (3), being nonresonant, is already a condition on the Taylor series of $\phi$.

A formal (holomorphic) transformation, sending $0$ to $0$, is a transformation of the form $(z',w')=(F(z,w),G(z,w))$, where $F$ and $G$ are formal power series in $z,w$ such that $F(0)=G(0)=0$. The formal mapping is {\it invertible} if its Jacobian determinant at $0$ is not zero. If the formal transformation is invertible, then we shall also refer to $(z',w')$ as {\it formal (holomorphic) coordinates} at $0$. We shall say that a formal transformation $(z',w')=(F(z,w),G(z,w))$ sends one formal hypersurface $M$ into another $M'$ if
\begin{equation}\Label{formalmap}
\Im G(z,u+i\phi)=\phi'(F(z,u+i\phi),\overline{F(z,u+i\phi)}, \Re G(z,u+i\phi)),
\end{equation}
where $\phi=\phi(z,\bar z,u)$, and $\phi$, $\phi'$ denote the formal graphing power series of $M$, $M'$, respectively. Our main result is the following:

\begin{Thm}\Label{main0} Let $M$ be a formal hypersurface through $0$ in $\bC^2$, satisfying Conditions $(1')$, $(2')$, and $(3)$. Then there are formal holomorphic coordinates $(z,w)$ at $0$ such that $M$ is given as a formal graph
\begin{equation}
\Im w=\phi(z,\bar z,\Re w),
\end{equation}
where the formal (Hermitian) power series $\phi(z,\bar z,u)$ is of the form
\begin{equation}
\phi(z,\bar z,u)=u\left(|z|^2+\sum_{a,b\geq 0}N_{ab}(u)z^a\bar z^b\right),\quad N_{ab}(u)=\overline{N_{ba}(u)},\ u\in \bR,
\end{equation}
satisfying the following normalization conditions
\begin{equation}\Label{normalform10}
N_{a0}(u)=N_{a1}(u)=0\quad a\geq 0,
\end{equation}
and
\begin{equation}\Label{normalform20}
\frac{d N_{22}}{d u}(u)=\frac{d N_{32}}{d u}(u)=\frac{d N_{33}}{d u}(u)=0.
\end{equation}
Moreover, the only invertible formal transformations
\begin{equation}\Label{transagain0}
z'=F(z,w),\quad w'=G(z,w)
\end{equation}
that preserves the normalization \eqref{normalform10} and \eqref{normalform20} are of the form
\begin{equation}
F(z,w)=\alpha z,\ G(z,w)=sw,\quad \alpha\in S^1,\ s\in \bR^*.
\end{equation}

More generally, if $M$ only satisfies $(1')$ and $(2')$,
we can still obtain \eqref{normalform10} as well as
the non-resonant part of \eqref{normalform20}, i.e.\
$$\frac{d^{k-1} N_{22}}{d u^{k-1}}(0)=\frac{d^{k-1} N_{32}}{d u^{k-1}}(0)=\frac{d^{k-1} N_{33}}{d u^{k-1}}(0)=0$$
for all non-resonant $k\geq 2$.
\end{Thm}

For a formal hypersurface $M$ through $0$ in $\bC^2$, we shall denote by $\Aut_0(M)$ the stability group of $M$ at $0$, i.e., the group of invertible formal transformations that preserve $M$ at $0$. An immediate consequence of Theorem \ref{main0} is the following:

\begin{Cor}\Label{maincor} Let $M$ be a formal hypersurface through $0$ in $\bC^2$, satisfying Conditions $(1')$, $(2')$, and $(3)$. Then, $\Aut_0(M)$ is a subgroup of $S^1\times \bR^*$.
\end{Cor}

The realization of $\Aut_0(M)$ as a subgroup of $S^1\times \bR^*$ goes, of course, via the correspondence $(F,G)\mapsto (F_z(0),G_w(0))$ in normal coordinates. There is a vast literature of investigations concerning $\Aut_0(M)$ for CR manifolds, but most treat $M$ only at finite type points. Papers that investigate $\Aut_0(M)$ at infinite type points include \cite{ELZ03}, \cite{Kow03}, \cite{Kow05}, \cite{ELZ09}, \cite{JuhlinLamel13}, \cite{KoL14}, \cite{KoL15}. The results in Corollary \ref{maincor} are more precise (for the class of manifolds considered here) than the results contained in these papers.

This paper is organized as follows. In the next section (which is broken into three subsections), the setup and normalization procedure is described and subsequently summarized in Theorem \ref{main1} (which readily translates into Theorem \ref{main0}). In the third and last section, some examples are given.

\section{Normalization}

\subsection {Setup}\Label{setup} Let $M$ be a formal hypersurface through $0$ in $\bC^2$. It is well known (see e.g.\ \cite{BER99a}) that there are formal holomorphic coordinates $(z,w)$ at $0$ such that $M$ is given by a graphing equation
\begin{equation}\Label{defeq}
\Im w=\phi(z,\bar z, \Re w),
\end{equation}
where $\phi(z,\chi,u)$ is a formal power series in $(z,\chi,u)$, which is Hermitian in $(z,\chi)$, i.e.
\begin{equation}\Label{Hermphi}
\phi(\bar\chi,\bar z,u)=\overline{\phi(z,\chi,u)}
\end{equation}
and further satisfies the normalization
\begin{equation}\Label{orgnorm}
\phi(z,0,u)=\phi(0,\chi,u)\equiv 0.
\end{equation}
We shall assume in this paper that $M$ is of infinite type at $0$ (i.e., Condition (1') above), which manifests itself in the defining equation \eqref{defeq} by $\phi(z,\bar z,u)$ satisfying
\begin{equation}\Label{inftype}
\phi(z,\bar z,0)\equiv 0.
\end{equation}
Equation \eqref{inftype} implies that the formal powers series $\phi(z,\bar z,u)$ has the following form
\begin{equation}\Label{phiexp}
\phi(z,\bar z,u)=\sum_{a,b\ge 0}\phi_{ab}z^{a}\bar z^{b}u + \sum_{a,b\ge 0,c\ge 2}\phi_{abc}z^{a}\bar z^{b} u^{c},
\end{equation}
and equation \eqref{orgnorm} implies that
\begin{equation}\Label{orgnorm2}
\phi_{0b}=\phi_{a0}=0,\quad \phi_{0bc}=\phi_{a0c}=0,\quad \forall a,b\geq 0.
\end{equation}
We shall consider the class of infinite type hypersurfaces that also satisfy Condition (2') above, which here is equivalent to $\phi_{11}\neq 0$. It is easy to see that a linear transformation in the $z$-variable will make $\phi_{11}=1$.

\subsection{Preliminary normalization}\Label{prenormsec} We shall normalize the defining equation of $M$ further by making formal transformations of the form
\begin{equation}\Label{trans}
z'=z+f(z,w),\quad w'=w+g(z,w),
\end{equation}
where $f(z,w)$ is a power series without constant term or linear term in $z$, and $g(z,w)$ a power series without constant term or linear term in $w$. Thus, we have the expansions
\begin{equation}\Label{expftrans}
f(z,w)=\sum_{l,k\geq0}f_{lk}z^lw^k,\quad f_{00}=f_{10}=0,
\end{equation}
and
\begin{equation}\Label{expgtrans}
g(z,w)=\sum_{l,k\geq0}g_{lk}z^lw^k,\quad g_{00}=g_{01}=0.
\end{equation}
We shall assume that $M$ is initially given in the $(z',w')$ coordinates at $0$ by
\begin{equation}\Label{defeq'}
\Im w'=\phi'(z',\bar z', \Re w'),
\end{equation}
with expansion
\begin{equation}\Label{phiexp'}
\phi'(z',\bar z',u')=\sum_{a,b\ge 0}\phi'_{ab}(z')^{a}(\bar z')^{b}u' + \sum_{a,b\ge 0,c\ge 2}\phi'_{abc}(z')^{a}(\bar z')^{b} (u')^{c},
\end{equation}
satisfying the prenormalization described above, i.e.,
\begin{equation}\Label{orgnorm2'}
\phi'_{11}=1,\quad \phi'_{0b}=\phi'_{a0}=0,\quad \phi'_{0bc}=\phi'_{a0c}=0,\quad \forall a,b\geq 0.
\end{equation}
We shall now describe how the transformation \eqref{trans} affects the coefficients in the defining equation. In the new coordinates $(z,w)$, the hypersurface $M$ is given by the equation \eqref{defeq} and we have the basic equation
\begin{equation}\Label{basic}
\phi+\Im g(z,u+i\phi)=
\phi'\big(z+f(z,u+i\phi),\bar z+\1{f(z,u+i\phi)}, u+ \Re g(z,u+i\phi)\big),
\end{equation}
where $\phi=\phi(z,\bar z,u)$. We shall only make transformations \eqref{trans} that preserve the prenormalization \eqref{orgnorm2'}. (Note that $\phi'_{11}=1$ is always preserved by the form of the transformation in \eqref{trans}.) It is well known (see \cite{BER99a}) that the prenormalization \eqref{orgnorm} holds in the coordinates $(z,w)$ if and only if a defining equation $\rho(z,\bar z,w,\bar w)=0$ of $M$ at $0$ satisfies $\rho(z,0,w,w)\equiv 0$; in our context, this means that
\begin{equation}\Label{normpres}
\frac{1}{2i}(g(z,w)-\bar g(0,w))=\phi'\left(z+f(z,w),\bar f(0,w), w+(g(z,w)-\bar g(0,w))/2\right)
\end{equation}
holds identically,
where the notation $\bar h(z,w)=\overline{h(\bar z,\bar w)}$ is used. We shall return to this characterization of this prenormalization in Lemma \ref{glkdetermine} below. For now, we just note some immediate conditions on $g(z,w)$ imposed by \eqref{normpres}. Setting $w=0$ in this identity, we see that $g(z,0)\equiv 0$, i.e.,
\begin{equation}\Label{gl0}
g_{l0}=0,\quad l\geq 0.
\end{equation}
Next, identifying coefficients of the monomial $z^lw$ in \eqref{normpres}, using $\phi'_{a0}=0$ and \eqref{gl0}, it is not difficult to see that we also have
\begin{equation}\Label{gl1}
g_{l1}=0,\quad l\geq 0,
\end{equation}
since every term in the expansion of the right hand side of \eqref{normpres} has at least two powers of $w$.

We also expand $\phi(z,\bar z,u)$ as in \eqref{phiexp}; the prenormalization implies that \eqref{orgnorm2} holds, and the form of the transformation \eqref{trans} guarantees that we retain the identity $\phi_{11}=1$. We shall now normalize $\phi$ further. We use the notation
\begin{equation}
\Delta \phi_{ab}:=\phi_{ab}-\phi'_{ab},\quad \Delta \phi_{abc}:=\phi'_{abc}-\phi_{abc},
\end{equation}
Our first lemma is the following, in which we assume that the prenormalization is preserved.

\begin{Lem}\Label{transf0}
For a fixed $l\geq 2$, the following transformation rule holds, modulo a non-constant term polynomial in $f_{a0}$, with $a<l$, whose coefficients depend only on the coefficients of $\phi'(z',\bar z',u')$:
\begin{equation}\Label{CMprenorm}
\Delta\phi_{l1}=-f_{l0}.
\end{equation}
\end{Lem}

\begin{proof} If we identify coefficients of $z^l\bar z u$ in the expansion of \eqref{basic}, then we get $\phi_{l1}$ only on the left hand side in view of \eqref{gl1}. Let us examine the right hand side. We note that $\phi'(z',\bar z',\bar u')$ has at least one power of $u'$, and $u+\Re g(z,u+i\phi)$ has at least one power of $u$ and any term involving $g$ has at least two powers of $u$ and cannot contribute to a term $z^l\bar z u$. A factor $\bar z$ can only come from $\bar z+\overline{f(z,u+i\phi)}$ and $\bar f$ will contribute another power of $u$. Thus, the only terms of the form $z^l\bar z u$ on the right will be from $\phi'_{a1}(z+f(z,0))^a\bar z u$ for $a\leq l$. Since $\phi'_{11}=1$, the conclusion of Lemma \ref{transf0} follows.
\end{proof}

It follows that we may perform an additional initial normalization of the defining equation of $M$ and require, in addition to the prenormalization above, that $\phi_{l1}=0$ for $l\geq 2$. In what follows, we shall assume that this is part of the prenormalization, i.e., in addition to \eqref{orgnorm2'}, we also assume
\begin{equation}\Label{addinnorm'}
\phi'_{l1}=0,\quad l\geq 2,
\end{equation}
and we shall consider only transformations that preserve this form. It follows from Lemma \ref{transf0} (and the fact that $f$ has no constant term or linear term in $z$) that we must impose
\begin{equation}\Label{fl0}
f_{l0}=0,\quad l\geq 0.
\end{equation}
It is not difficult to see that if we require \eqref{fl0}, then
\begin{equation}\Label{fl0'}
\Delta \phi_{ab}=0
\end{equation}
for all $a,b$. For convenience of notation, we shall therefore drop the $'$ on $\phi'_{ab}$ and simply write $\phi_{ab}$.

We also have the following lemma:

\begin{Lem}\Label{flkdetermine} For fixed $l\geq 3$, $k\geq 2$, the following transformation rule holds, modulo a non-constant term polynomial in $f_{a,b-1},\bar f_{a,b-1}, g_{ab}, \bar g_{0b},\phi_{a1b},\bar \phi_{a1b}$, with $b<k$, whose coefficients depend on the coefficients of $\phi'(z',\bar z',u')$:
\begin{equation}\Label{CMprenorm2}
\Delta\phi_{l1k}=\frac{k-1}{2}\,g_{l-1,k}-f_{l,k-1}-2\,\phi'_{l2}\bar f_{0,k-1}.
\end{equation}
\end{Lem}

\begin{proof} We identify the coefficients of $z^l\bar z u^k$ in the expansion of \eqref{basic}. By examining the expansion of $\Im g(z,u+i\phi)$ and using the prenormalization conditions, we observe that on the left hand side we get
\begin{equation}\Label{LHSflk}
\phi_{l1k}+\frac{k}{2}\,g_{l-1,k}+\sum_
{2\leq c'\leq k-1\atop 1\leq a'\leq l-1}
\frac{k+1-c'}{2}\,g_{l-a',k+1-c'}\phi_{a'1c'}+\sum_{2\leq c'\leq k-1}(k+1-c')\,\Re g_{0,k+1-c'}\phi_{l1c'},
\end{equation}
which is equal to
$$
\phi_{l1k}+\frac{k}{2}g_{l-1,k}
$$
modulo a non-constant term polynomial in $g_{ab},\bar g_{0b}, \phi_{a'1c'}, \bar \phi_{a'1c'}$, with $a, a'< l$, $b,c'<k$.
On the right hand side, we examine the term (with the understanding that $\phi'_{ab1}=\phi'_{ab}=\phi_{ab}$)
\begin{equation}\Label{phi'abc}
\phi'_{abc}(z+f(z+i\phi))^a(\bar z+\overline{f(z,u+i\phi)})^b(u+\Re g(z,u+i\phi))^c
\end{equation}
and first observe that if any term from the expansion of $\phi(z,\bar z,u)$ is involved, then it can only be of the form $\phi_{a1b}$, which contributes one power of $z$ and $b$ powers of $u$. The contribution from $\bar z+\overline{f(z,u+i\phi)}$ can then only be through a factor of the form $\bar f_{0b'}$, which will contribute at least one power of $u$. Since the term $u+\Re g(z,u+i\phi)$ always contributes at least a factor of $u$ as well, we conclude (after some thought) that the term that involves $\phi$ will be a non-constant term polynomial in $f_{a,b-1},\bar f_{0,b-1}, g_{ab}, \bar g_{0b},\phi_{a1b},\bar \phi_{a1b}$, with $b<k$. If $\phi$ is not involved in the term on the right, then we can only get the single power $\bar z$ from $\bar z+\overline{f(z,u+i\phi)}$. If $c\geq 2$ in \eqref{phi'abc}, then we of course get the term $\phi'_{l1k}$ if $c=k$ and we pick the single term that does not involve $f$, $\bar f$, or $\Re g$. If the latter are involved, we note that $u+\Re g(z,u+i\phi)$ contributes at least two powers of $u$ and no term from $\Re g(z,u+i\phi)$ can involve a $g_{ab}$ with $b\geq k-1$. We conclude that any contribution from $f$ or $\bar f$ must be through a factor of $f_{a,b-1},\bar f_{a,b-1}$ with $b<k$. Thus, beside the term $\phi'_{l1k}$, the terms arising from \eqref{phi'abc}, with $c\geq 2$, will be a non-constant term polynomial in $f_{a,b-1},\bar f_{0,b-1}, g_{ab}, \bar g_{0b}$, with $b<k$. Finally, if $c=1$ in \eqref{phi'abc}, then we get the contribution
$$
f_{l,k-1}+\frac{1}{2}g_{l-1,k}
$$
from the term with $a=b=1$ (recall $\phi_{11}=1$ and $\phi'_{ab}=\phi_{ab}$ by our prenormalization) and the contribution $2\phi_{l2}\bar f_{0,k-1}$ from the term with $a=l,b=2$, but the remaining terms will be a nonconstant term polynomial in $f_{a,b-1},\bar f_{a,b-1}, g_{ab}, \bar g_{ab}$, with $b<k$.
 This completes the proof.
\end{proof}

We now return to see what the characterization \eqref{normpres} of the prenormalization \eqref{orgnorm2'} yields for the coefficients $g_{lk}$:
\begin{Lem}\Label{glkdetermine} For each $k\geq 2$ there are
\begin{itemize}
\item[(i)] a non-constant term polynomial $P_k$ in the variables $f_{0,b-1}$, $\bar f_{0,b-1}$, and $\Re g_{0b}$, with $b<k$, whose coefficients depend only on the coefficients of $\phi'(z',\bar z',u')$, and

\item[(ii)] for each $l\geq 1$, a non-constant term polynomial $Q_{lk}$ in $f_{a,b-1}$, $\bar f_{0,b-1}$, $g_{ab}$, and $\bar g_{0b}$, with $a\leq l$ and $b<k$, whose coefficients depend only on the coefficients of $\phi'(z',\bar z',u')$,
\end{itemize}
with the following property: The transformation \eqref{trans} preserves the prenormalization \eqref{orgnorm2'} if and only if $\Im g_{0k}=0$ modulo the value of $P_k$ for every $k\geq 2$,  and $g_{1k}=\bar f_{0,k-1}$, $g_{lk}=0$ both modulo the value of $Q_{lk}$ for every $l\geq 2$ and $k\geq 2$.
\end{Lem}

\begin{proof} To find $P_k$ in (i), we identify coefficients of the monomial $w^k$ in \eqref{normpres}. On the left hand side, we get $\Im g_{0k}$. On the right hand side, we note that $\phi'(z',\bar z',u')$ contributes at least one power each of $z',\bar z',u'$. It is clear that any term in the expansion of the right hand side of \eqref{normpres} that contributes $w^k$ will have a coefficient that is a product of $f_{0,b-1}$, $\bar f_{0,b-1}$, and $\Re g_{0b}$, with $b<k$, and a coefficient from the expansion of $\phi'(z',\bar z',u')$. This establishes the existence of $P_k$ in (i).

To find $Q_{lk}$ in (ii), we identify coefficients of the monomial $z^lw^k$ in \eqref{normpres}. The only contribution on the left hand side is $g_{lk}/2i$. Since we are also requiring the normalization \eqref{addinnorm'} and have already established \eqref{gl0}, \eqref{gl1}, \eqref{fl0}, every contribution to the coefficient of $z^lw^k$ is a product of $f_{a,b-1}$, $\bar f_{0,b-1}$, $g_{ab}$, $\bar g_{0b}$, with $b<k$, and a coefficient in the expansion of $\phi'(z',\bar z',u')$, {\it} except when $l=1$ in which case there is a coefficient of the form $\bar f_{0,k-1}$; for $l\geq 2$, the analogous term
$
\phi_{l1}\bar f_{0,k-1}
$
vanishes by \eqref{addinnorm'}.

The conclusion in (ii) now follows.
\end{proof}

\begin{Rem}{\rm The conditions on the coefficients $g_{lk}$ in Lemma \ref{glkdetermine} above can also be derived by considering the transformation rules for $\Delta_{l0k}$ stemming from \eqref{basic}. For reasons that will become apparent in the next section, it will be convenient to do so specifically for the coefficients $g_{1k}$.
}
\end{Rem}

Lemma \ref{flkdetermine} suggests an additional prenormalization ($\phi'_{l1k}=0$, $l\geq 3$, $k\geq 2$), and this lemma together with Lemma \ref{glkdetermine} leads to an induction scheme that can be summarized in the following proposition:

\begin{Pro}\Label{transcoeffs} In addition to the prenormalization given by \eqref{orgnorm2'} and \eqref{addinnorm'}, the following prenormalization
\begin{equation}\Label{addinnorm'2}
\phi'_{l1k}=0,\quad l\geq 3,\ k\geq 2,
\end{equation}
can be achieved. Any transformation of the form \eqref{trans} that preserves the prenormalization given by \eqref{orgnorm2'}, \eqref{addinnorm'}, and \eqref{addinnorm'2}, satisfies \eqref{gl0}, \eqref{gl1}, \eqref{fl0}, and the coefficients
\begin{equation}\Label{transcoeffe1}
f_{l+1,k-1},\ \Im g_{0k},\ g_{lk},\quad k,l\geq 2,
\end{equation}
are given by non-constant term polynomials, whose coefficients depend only on the coefficients of $\phi'(z',\bar z', u')$, in the variables
\begin{equation}\Label{transvars}
f_{0,k-1},\ f_{1,k-1},\ f_{2,k-1},\ \Re g_{0k},\ g_{1k},\quad k\geq 2,
\end{equation}
and their complex conjugates.
\end{Pro}

\begin{proof} The proof is a straightforward induction on $k\geq 2$ using Lemmas \ref{glkdetermine} and \ref{flkdetermine}, together with the normalizations \eqref{gl0}, \eqref{gl1}, \eqref{fl0}. Details are left to the reader.
\end{proof}

\subsection{Complete normalization}\Label{completenormsec} Our aim now is to find a final (complete) normalization of the defining equation of $M$ at $0$  that uniquely determines the variables in \eqref{transvars}, and has the property that this normalization is preserved only when these variables vanish. Proposition \ref{transcoeffs} will then imply that the only transformation of the form \eqref{trans} that preserves this complete normalization is the identity mapping.

We shall assume that the prenormalizations in the previous subsection are preserved (although as alluded to in that section, we will also study the transformation rules for $\Delta\phi_{10k}$). Recall that our prenormalization implies that $\phi'_{ab}=\phi_{ab}$ and we will drop $'$ on $\phi'_{ab}$ so simplify the notation. The main technical result in this paper is contained in the following lemma.

\bl\Label{transf}
The transformation rules are given as follows, modulo non-constant term
polynomials, whose coefficients are given by the coefficients in the expansion of $\phi'(z',\bar z',u')$, in the variables consisting of $f_{a,b-1}, g_{ab}$ and their complex conjugates, with $b<k$:
\begin{eqnarray*}
\D\phi_{10k} =& \frac1{2i} g_{1k} - \1f_{0,k-1}\\
\D\phi_{11k}  =& (k-1)\Re g_{0k} - 2\Re f_{1,k-1}\\
\D\phi_{21k} =&  (\frac{k}2-\frac12) g_{1k}
+(i(k-1)-2\phi_{22}) \bar f_{0,k-1} - f_{2,k-1}\\
\D\phi_{22k} =&  -6\Re( \phi_{23} \bar f_{0,k-1})\\
&+ (k-1)\phi_{22} \Re g_{0k} +2(k-1)\Im f_{1,k-1} -4\phi_{22}\Re f_{1,k-1} ,\\
\D\phi_{32k} =& (\frac{k-1}2\phi_{22}+ \frac{i}2{k\choose2} - \frac{ik}2)g_{1k}
+({k-1\choose2}+3i(k-1)\phi_{22} - 3\phi_{33})\bar f_{0,k-1}   - 4\phi_{42}f_{0,k-1}\\
 &+(k-1) \phi_{32}\Re g_{0k}
 - \phi_{32}(5\Re f_{1,k-1}+\Im f_{1,k-1}) -(i(k-1)+\phi_{22})f_{2,k-1}\\
\D\phi_{33k} =& \Re((k-1)\phi_{23} g_{1k}) -  8\Re (((k-1)i\phi_{23} - \phi_{34})\bar f_{0,k-1} )
\\&+((k-1)\phi_{33}-{k\choose3} +{k\choose2})\Re g_{0k} +
2({k-1\choose2} -3\phi_{33} )\Re f_{1,k-1} +  6(k-1)\phi_{22} \Im f_{1,k-1}\\
& - 4\Re(\phi_{23} f_{2,k-1}),
\end{eqnarray*}
where $k\ge 2$ and we use the convention ${a\choose b}=0$ whenever $a<b$.
\el

\bpf
We begin by inspecting the terms arising from the expansion of the left-hand side of \eqref{basic}:
\begin{equation}
\sum \Im g_{lk}z^{l}(u+i\sum\phi_{abc})^{k}.
\end{equation}
The only term with no $\phi_{abc}$ that is relevant to the transformation rules in Lemma~\ref{transf}
is $\Im g_{1k}zu^{k}$ which contributes to $\D\phi_{10k}$
with coefficient $\frac1{2i}$. Furthermore, factors $\phi_{abc}$ with $c\ge2$
cannot contribute since they make the total degree in $u$ greater than $k$.
It remains to consider terms with one or several
factors $\phi_{ab}$. In view of the preservation of \eqref{orgnorm2'} and \eqref{addinnorm'}, these
can only be
\begin{equation}
\phi_{11}=1,\quad \phi_{22},\quad\phi_{23},\quad\phi_{32},\quad\phi_{33}.
\end{equation}
Recall that, by Lemma \ref{glkdetermine}, all
coefficients of $g$ except $\Re g_{0k}$ and $g_{1k}$
are determined by $f_{a,b-1}, g_{ab}$, and their complex conjugates, with $b<k$.

We next inspect terms with $g_{1k}$ that appear as
\begin{equation}
\frac1{2i}g_{1k}z(u+i\sum\phi_{ab}z^{a}\bar z^{b}u)^{k}.
\end{equation}
The term with single factor $\phi_{11}$ contributes
as $\frac1{2i}ik  g_{1k} z^{2}\bar z u^{k}$ to $\D\phi_{21k}$.
The term with single factor $\phi_{22}$ contributes
as $\frac1{2i}ik\phi_{22}  g_{1k} z^{3}\bar z^{2} u^{k}$ to $\D\phi_{32k}$.
The term with single factor $\phi_{23}$ contributes
as $\frac1{2i}ik\phi_{23}  g_{1k} z^{3}\bar z^{3} u^{k}$ to $\D\phi_{33k}$.
The term with single factor $\phi_{32}$ contributes
as its conjugate $\frac1{2i}ik\phi_{32}  \bar g_{1k} z^{3}\bar z^{3} u^{k}$ to $\D\phi_{33k}$.
The factor $\phi_{33}$ has no contribution to the identities in the lemma.
Further, the term with the square of $\phi_{11}$ contributes
as $\frac1{2i}i^{2}{k\choose2}  g_{1k} z^{3}\bar z^{2} u^{k}$ to $\D\phi_{32k}$.
Other products $\phi_{ab}\phi_{cd}$ have no contribution.
Also terms with more than $2$ factors $\phi_{ab}$ have no contribution.

Next consider terms with $\Re g_{0k}$ that appear as
\begin{equation}
\Im (u+i\sum\phi_{ab}z^{a}\bar z^{b}u)^{k} \Re g_{0k}.
\end{equation}
The term with single factor $\phi_{11}$ contributes
as $k z\bar z u^{k} \Re g_{0k}$ to $\D\phi_{11k}$.
The term with single factor $\phi_{22}$ contributes
as $k\phi_{22} z^{2}\bar z^{2} u^{k} \Re g_{0k}$ to $\D\phi_{22k}$.
The terms with single factors $\phi_{32}$ and $\phi_{23}$ contribute
as $\Im(ik\phi_{32} z^{3}\bar z^{2}+ik\phi_{23}z^{2}\bar z^{3}) u^{k} \Re g_{0k}$ to $\D\phi_{32k}$.
The term with single factor $\phi_{33}$ contributes
as $k\phi_{33} z^{3}\bar z^{3} u^{k} \Re g_{0k}$ to $\D\phi_{33k}$.
Next, there is no contribution from terms with products $\phi_{ab}\phi_{cd}$
because of the reality of $\phi$.
Finally, the term with the cube of $\phi_{11}$ contributes
as $\Im (i^{3}{k\choose3} z^{3}\bar z^{3} u^{k}) \Re g_{0k}$ to $\D\phi_{33k}$.
Other terms have no contribution.

We now inspect the terms on the right-hand side of \eqref{basic} that contribute with minus.
Those containing $f_{l,k-1}$ and $g_{lk}$ arise from the expansion of
\begin{equation}
-\sum\phi_{ab}(z+\sum f_{l,k-1}z^{l}(u+i\phi)^{k-1})^{a}(\bar z+\sum \bar f_{l,k-1}\bar z^{l}(u-i\phi)^{k-1})^{b}
(u+\Re \sum g_{lk}z^{l}(u+i\phi)^{k}),
\end{equation}
where the $\phi_{ab}$ that occur are technically $\phi'_{ab}$ but we recall that $\phi'_{ab}=\phi_{ab}$ as a consequence of \eqref{fl0}.

We first collect the terms with $g_{1k}$ that appear as
\begin{equation}
-\phi_{ab}z^{a}\bar z^{b} \frac12 z(u+i\phi)^{k}  g_{1k}.
\end{equation}
For $(a,b)=(1,1)$ we obtain $-\frac12 g_{1k}$ contributing to $\D\phi_{21k}$ and
$-\frac12 ik\phi_{11} g_{1k}=-\frac{ik}2 g_{1k}$ contributing to $\D\phi_{32k}$.
For $(a,b)=(2,2)$ we obtain $-\phi_{22}\frac12 g_{1k}$ contributing to $\D\phi_{32k}$.
For $(a,b)=(2,3)$ we obtain $-\phi_{23}\frac12 g_{1k}$ contributing to $\D\phi_{33k}$.
For $(a,b)=(3,2)$ we obtain its conjugate $-\phi_{32}\frac12 \bar g_{1k}$ contributing to the same term.
Other terms have no contribution.

We next consider the terms with $g_{0k}$ that appear as
\begin{equation}
-\phi_{ab}z^{a}\bar z^{b} \Re (u+i\phi)^{k}  \Re g_{0k}.
\end{equation}
For $(a,b)=(1,1)$ we obtain $-\Re g_{0k}$ contributing to $\D\phi_{11k}$ and
${k\choose2}\phi_{11}^{2} g_{0k}={k\choose2} g_{0k}$ contributing to $\D\phi_{33k}$.
For $(a,b)=(2,2)$ we obtain $-\phi_{22}\Re g_{0k}$ contributing to $\D\phi_{22k}$.
For $(a,b)=(3,2)$ we obtain $-\phi_{32}\Re g_{0k}$ contributing to $\D\phi_{32k}$.
For $(a,b)=(3,3)$ we obtain $-\phi_{33}\Re g_{0k}$ contributing to $\D\phi_{33k}$.

As our final consideration we deal with the terms involving $f_{l,k-1}$.
We begin with terms involving $f_{0,k-1}$ that arise as
\begin{equation}
-\phi_{ab} (a z^{a-1}\bar z^{b}f_{0,k-1}(u+i\phi)^{k-1} + b z^{a}\bar z^{b-1}\bar f_{0,k-1}(u-i\phi)^{k-1})u.
\end{equation}

For $(a,b)=(1,1)$ we obtain $-\bar f_{0,k-1}$ contributing to $\D\phi_{10k}$,
 $-\phi_{11}\bar f_{0,k-1}(-i)(k-1)\phi_{11}=i(k-1)\bar f_{0,k-1}$ contributing to $\D\phi_{21k}$,
 $-\phi_{11}\bar f_{0,k-1}(-i)^{2} {k-1\choose 2}\phi_{11}^{2}={k-1\choose2}\bar f_{0,k-1}$
 and $-\phi_{11}\bar f_{0,k-1}(-i)(k-1)\phi_{22}=i(k-1)\phi_{22}\bar f_{0,k-1}$ both
  contributing to $\D\phi_{32k}$, and
  $-\phi_{11}\bar f_{0,k-1}(-i)(k-1)\phi_{23}=i(k-1)\phi_{23} \bar f_{0,k-1}$ and its conjugate
   $-\phi_{11} f_{0,k-1}i(k-1)\phi_{32}=-i(k-1)\phi_{32} f_{0,k-1}$  both contributing to
   $\D\phi_{33k}$.

Next, for $(a,b)=(2,2)$ we obtain
$-2\phi_{22} \bar f_{0,k-1}$ contributing to $\D\phi_{21k}$,
$-2\phi_{22} (k-1)(-i)\phi_{11} \bar f_{0,k-1}=2i(k-1)\phi_{22} \bar f_{0,k-1}$
contributing to $\D\phi_{32k}$.

For $(a,b)=(2,3)$ and $(a,b)=(3,2)$ we obtain
$-3\phi_{23}\bar f_{0,k-1}$ and its conjugate
$-3\phi_{32} f_{0,k-1}$ contributing to $\D\phi_{22k}$,
$-3\phi_{23}\bar f_{0,k-1}(k-1)(-i\phi_{11})=3(k-1)i\phi_{23}\bar f_{0,k-1}$
and its conjugate
$-3\phi_{32} f_{0,k-1}(k-1)(i\phi_{11})=-3(k-1)i\phi_{32}f_{0,k-1}$
contributing to $\D\phi_{33k}$.

For $(a,b)=(3,3)$ we obtain
$-3\phi_{33}\bar f_{0,k-1}$ contributing to $\D\phi_{32k}$.

For $(a,b)=(4,2)$ we obtain
$-4\phi_{42} f_{0,k-1}$ also contributing to $\D\phi_{32k}$.

For $(a,b)=(3,4)$ and $(a,b)=(4,3)$ we obtain
$-4\phi_{34}\bar f_{0,k-1}$ and its conjugate
$-4\phi_{43} f_{0,k-1}$ both contributing to $\D\phi_{33k}$.

Other terms have no contribution.

We next treat terms involving $f_{1,k-1}$ that arise as
\begin{equation}
-\phi_{ab} (af_{1,k-1}(u+i\phi)^{k-1} + b\bar f_{1,k-1}(u-i\phi)^{k-1}) z^{a}\bar z^{b}u.
\end{equation}

For $(a,b)=(1,1)$ we obtain
$-f_{1,k-1}-\bar f_{1,k-1}=-2\Re f_{1,k-1}$
contributing to $\D\phi_{11k}$,
$$-\phi_{11} (f_{1,k-1}(k-1)i\phi_{11} + \bar f_{1,k-1}(k-1)(-i\phi_{11}))
=2(k-1)\Im f_{1,k-1}$$ contributing to $\D\phi_{22k}$,
and
$$-\phi_{11} (f_{1,k-1}(k-1)i\phi_{22} + \bar f_{1,k-1}(k-1)(-i\phi_{22}))
=-2(k-1)\Re(i\phi_{22} f_{1,k-1})$$
and
$$-\phi_{11} (f_{1,k-1}{k-1\choose2}(i\phi_{11})^2 + \bar f_{1,k-1}{k-1\choose2}(-i\phi_{11})^2)
=2{k-1\choose2}\Re f_{1,k-1}$$
 both contributing to $\D\phi_{33k}$.

For $(a,b)=(2,2)$ we obtain
$-\phi_{22}(2f_{1,k-1}+2\bar f_{1,k-1})=-4\phi_{22}\Re f_{1,k-1}$
contributing to $\D\phi_{22k}$ and
$$-\phi_{22} (2f_{1,k-1}(k-1)i\phi_{11} + 2\bar f_{1,k-1}(k-1)(-i\phi_{11}))
=4(k-1)\phi_{22}\Im f_{1,k-1}$$ contributing to $\D\phi_{33k}$.

For $(a,b)=(3,2)$ we obtain
$-\phi_{32}(3f_{1,k-1}+2\bar f_{1,k-1})=-\phi_{32}(5\Re f_{1,k-1} + \Im f_{1,k-1})$
contributing to $\D\phi_{32k}$.

Finally for $(a,b)=(3,3)$ we obtain
$-\phi_{32}(3f_{1,k-1}+3\bar f_{1,k-1})=-6\phi_{33}\Re f_{1,k-1}$
contributing to $\D\phi_{33k}$.

It remains to deal with terms involving $f_{2,k-1}$ that arise as
\begin{equation}
-\phi_{ab} (af_{2,k-1}z(u+i\phi)^{k-1} + b\bar f_{2,k-1}\bar z(u-i\phi)^{k-1}) z^{a}\bar z^{b}u.
\end{equation}

For $(a,b)=(1,1)$ we obtain
$-f_{2,k-1}$
contributing to $\D\phi_{21k}$, and
$$-f_{2,k-1}(k-1)(i\phi_{11})=-i(k-1)f_{2,k-1}$$
contributing to $\D\phi_{32k}$.

For $(a,b)=(2,2)$ we obtain
$-\phi_{22}f_{2,k-1}$
contributing to $\D\phi_{32k}$.

For $(a,b)=(2,3)$ we obtain
$-\phi_{23} 2f_{2,k-1}$
and for $(a,b)=(3,2)$ its conjugate
$ - \phi_{32} 2\bar f_{2,k-1}$
both contributing to $\D\phi_{33k}$.
Other terms have no contribution.
\epf

Extracting real and imaginary parts we obtain from Lemma \ref{transf} the following identity, modulo a vector of non-constant term
polynomials, whose coefficients are given by the coefficients in the expansion of $\phi'(z',\bar z',u')$, in the variables consisting of $f_{a,b-1}, g_{ab}$ and their complex conjugates, with $b<k$:

\begin{equation}\Label{reimmatrix}
\begin{pmatrix}
\Re\D\phi_{10k}\\
\Im\D\phi_{10k}\\
\D\phi_{11k} \\
\Re\D\phi_{21k}\\
\Im\D\phi_{21k}\\
\D\phi_{22k}\\
\Re\D\phi_{32k}\\
\Im\D\phi_{32k}\\
\D\phi_{33k}
\end{pmatrix}
= A
\begin{pmatrix}
\Re g_{1k} \\
\Im g_{1k}\\
\Re f_{0,k-1}\\
\Im f_{0,k-1}\\
\Re g_{0k}\\
\Re f_{1,k-1}\\
\Im f_{1,k-1}\\
\Re f_{2,k-1}\\
\Im f_{2,k-1}
\end{pmatrix},
\end{equation}
where the matrix $A$ is

\footnotesize
$$
\begin{pmatrix}
0		& \frac1{2}& -1 	&0 		& 0	&0		&0&0&0\\
-\frac12	&0		&0	&1 		&0	&0		&0&0&0\\
0		&0		&0    &0 		&k-1	&-2		&0&0&0\\
\frac{k-1}2&0		&-2\phi_{22}&k-1		&0&0&0& -1&0\\
0		&\frac{k-1}2&k-1&2\phi_{22}	&0&0&0&0&-1\\
0		&0		&-6\Re\phi_{32}&6\Im\phi_{32}	&(k-1)\phi_{22}&-4\phi_{22}&2(k-1)&0&0\\
\frac{k-1}2\phi_{22}&-\frac{k(k-3)}4&
\begin{matrix}
\frac{k(k-1)}2
-3\phi_{33}\\
-4\Re\phi_{42}
\end{matrix}
&
3(k-1)\phi_{22} & (k-1)\Re\phi_{32}
&-5\Re\phi_{32}&-\Re\phi_{32}
&-\phi_{22}&k-1\\
\frac{k(k-3)}4&\frac{k-1}2\phi_{22} &3(k-1)\phi_{22}&
\begin{matrix}
-\frac{k(k-1)}2+3\phi_{33}\\
+4\Im\phi_{42}
\end{matrix}
&(k-1)\Im\phi_{32}
&-5\Im\phi_{32}&-\Im\phi_{32}
&-k+1&-\phi_{22}\\

\begin{matrix}
(k-1)\\
\times\Re\phi_{32}
\end{matrix}
&
\begin{matrix}
(k-1)\\
\times\Im\phi_{32}
\end{matrix}
&
\begin{matrix}
 8\Re\phi_{43}\\
 -8(k-1)\Im \phi_{32}
\end{matrix}
&
\begin{matrix}
 8\Im\phi_{43}\\
 +8(k-1)\Re \phi_{32}
\end{matrix}
&
\begin{matrix}
\frac{k(k-1)(6-k)}6\\
+(k-1)\phi_{33}
\end{matrix}
&
\begin{matrix}
k(k-1)\\-3\phi_{33}
\end{matrix}
&
\begin{matrix}
6(k-1)\\
\times\phi_{22}
\end{matrix}
&
-4\Re\phi_{32}
&-4\Im\phi_{32}
\end{pmatrix}.
$$
\normalsize
We have the following lemma:

\begin{Lem}\Label{detA} The determinant of $A$ is of the form
\begin{equation}\Label{detAdetB}
\det A=\frac{1}{4}(k-1)\det B,
\end{equation}
where $\det B$ is a polynomial in $k$ of degree $7$, whose leading coefficient is nonzero. Moreover, the coefficients of the polynomial $\det B$ depend only on $\phi_{ab}$ with $a,b\leq 4$ and $a+b\leq 7$.
\end{Lem}

\begin{Def}\Label{generalpos}
{\rm For a formal hypersurface $M\subset \bC^2$,
given by \eqref{defeq'} at $0\in M$ in any coordinate system satisfying the prenormalization described in Proposition \ref{transcoeffs}, we
define its {\em characteristic polynomial} $P(k, j^7_0\phi)$
to be $c\det B$, where $c$ is chosen such that $P$ is monic in $k$.
We call an integer $k\ge2$ a {\em resonance} for $M$ (at $0$) if $P(k,  j^7_0\phi)=0$.
Then $M$ is said to be {\it nonresonant}
if there are no resonances.
}
\end{Def}

\begin{proof}[Proof of Lemma $\ref{detA}$]
We shall apply elementary row operations to $A$.
Multiplying the first row by suitable factors and subtracting
from other rows below we obtain
\footnotesize
$$
\begin{pmatrix}
0		& \frac1{2}& -1 	&0 		& 0	&0		&0&0&0\\
-\frac12	&0		&0	&1 		&0	&0		&0&0&0\\
0		&0		&0    &0 		&k-1	&-2		&0&0&0\\
\frac{k-1}2&0		&-2\phi_{22}&k-1		&0&0&0& -1&0\\
0		&0		&2(k-1)&2\phi_{22}	&0&0&0&0&-1\\
0		&0		&-6\Re\phi_{32}&6\Im\phi_{32}	&(k-1)\phi_{22}&-4\phi_{22}&2(k-1)&0&0\\
\frac{k-1}2\phi_{22}&0&
\begin{matrix}
k
-3\phi_{33}\\

-4\Re\phi_{42}
\end{matrix}
&
3(k-1)\phi_{22} & (k-1)\Re\phi_{32}
&-5\Re\phi_{32}&-\Re\phi_{32}
&-\phi_{22}&k-1\\

\frac{k(k-3)}4&0 &4(k-1)\phi_{22}&
\begin{matrix}
-\frac{k(k-1)}2+3\phi_{33}\\
+4\Im\phi_{42}
\end{matrix}
&(k-1)\Im\phi_{32}
&-5\Im\phi_{32}&-\Im\phi_{32}
&-k+1&-\phi_{22}\\

\begin{matrix}
(k-1)\\
\times\Re\phi_{32}
\end{matrix}
&
0
&
\begin{matrix}
 8\Re\phi_{43}\\
 -6(k-1)\Im \phi_{32}
\end{matrix}
&
\begin{matrix}
 8\Im\phi_{43}\\
 +8(k-1)\Re \phi_{32}
\end{matrix}
&
\begin{matrix}
\frac{k(k-1)(6-k)}6\\
+(k-1)\phi_{33}
\end{matrix}
&
\begin{matrix}
k(k-1)\\-3\phi_{33}
\end{matrix}
&
\begin{matrix}
6(k-1)\\
\times\phi_{22}
\end{matrix}
&
-4\Re\phi_{32}
&-4\Im\phi_{32}
\end{pmatrix}.
$$
\normalsize

Similarly, multiplying the second row by suitable factors and subtracting
from other rows below we obtain
\footnotesize
$$
\begin{pmatrix}
0		& \frac1{2}& -1 	&0 		& 0	&0		&0&0&0\\
-\frac12	&0		&0	&1 		&0	&0		&0&0&0\\
0		&0		&0    &0 		&k-1	&-2		&0&0&0\\
0		&0		&-2\phi_{22}&2(k-1)		&0&0&0& -1&0\\
0		&0		&2(k-1)&2\phi_{22}	&0&0&0&0&-1\\
0		&0		&-6\Re\phi_{32}&6\Im\phi_{32}	&(k-1)\phi_{22}&-4\phi_{22}&2(k-1)&0&0\\

0		&0&
\begin{matrix}
k
-3\phi_{33}\\

-4\Re\phi_{42}
\end{matrix}
&
4(k-1)\phi_{22} & (k-1)\Re\phi_{32}
&-5\Re\phi_{32}&-\Re\phi_{32}
&-\phi_{22}&k-1\\

0	&0 &4(k-1)\phi_{22}&
\begin{matrix}
-k+3\phi_{33}\\
+4\Im\phi_{42}
\end{matrix}
&(k-1)\Im\phi_{32}
&-5\Im\phi_{32}&-\Im\phi_{32}
&-k+1&-\phi_{22}\\

0
&
0
&
\begin{matrix}
 8\Re\phi_{43}\\
 -6(k-1)\Im \phi_{32}
\end{matrix}
&
\begin{matrix}
 8\Im\phi_{43}\\
 +10(k-1)\Re \phi_{32}
\end{matrix}
&
\begin{matrix}
\frac{k(k-1)(6-k)}6\\
+(k-1)\phi_{33}
\end{matrix}
&
\begin{matrix}
k(k-1)\\-3\phi_{33}
\end{matrix}
&
\begin{matrix}
6(k-1)\\
\times\phi_{22}
\end{matrix}
&
-4\Re\phi_{32}
&-4\Im\phi_{32}
\end{pmatrix}.
$$
\normalsize

Now, multiplying the third row by suitable factors and subtracting
from other rows below we obtain
\footnotesize
$$
\begin{pmatrix}
0		& \frac1{2}& -1 	&0 		& 0	&0		&0&0&0\\
-\frac12	&0		&0	&1 		&0	&0		&0&0&0\\
0		&0		&0    &0 		&k-1	&-2		&0&0&0\\
0		&0		&-2\phi_{22}&2(k-1)		&0&0&0& -1&0\\
0		&0		&2(k-1)&2\phi_{22}	&0&0&0&0&-1\\
0		&0		&-6\Re\phi_{32}&6\Im\phi_{32}	&0&-2\phi_{22}&2(k-1)&0&0\\

0		&0&
\begin{matrix}
k
-3\phi_{33}\\

-4\Re\phi_{42}
\end{matrix}
&
4(k-1)\phi_{22} & 0
&-3\Re\phi_{32}&-\Re\phi_{32}
&-\phi_{22}&k-1\\

0&0 &4(k-1)\phi_{22}&
\begin{matrix}
-k+3\phi_{33}\\
+4\Im\phi_{42}
\end{matrix}
&0
&-3\Im\phi_{32}&-\Im\phi_{32}
&-k+1&-\phi_{22}\\

0
&
0
&
\begin{matrix}
 8\Re\phi_{43}\\
 -6(k-1)\Im \phi_{32}
\end{matrix}
&
\begin{matrix}
 8\Im\phi_{43}\\
 +10(k-1)\Re \phi_{32}
\end{matrix}
&
0
&
\begin{matrix}
\frac{k(2k+3)}3\\-\phi_{33}
\end{matrix}
&
\begin{matrix}
6(k-1)\\
\times\phi_{22}
\end{matrix}
&
-4\Re\phi_{32}
&-4\Im\phi_{32}
\end{pmatrix}.
$$
\normalsize

Next, repeating the procedure for the $4$th and $5$ row
to eliminate the $8$th and $9$th column respectively
we obtain
\footnotesize
$$
\begin{pmatrix}
0		& \frac1{2}& -1 	&0 		& 0	&0		&0&0&0\\
-\frac12	&0		&0	&1 		&0	&0		&0&0&0\\
0		&0		&0    &0 		&k-1	&-2		&0&0&0\\
0		&0		&-2\phi_{22}&2(k-1)		&0&0&0& -1&0\\
0		&0		&2(k-1)&2\phi_{22}	&0&0&0&0&-1\\
0		&0		&-6\Re\phi_{32}&6\Im\phi_{32}	&0&-2\phi_{22}&2(k-1)&0&0\\

0		&0&
\begin{matrix}
2k^{2}-3k+2+2\phi_{22}^{2}\\
-3\phi_{33}
-4\Re\phi_{42}
\end{matrix}
&
4(k-1)\phi_{22} & 0
&-3\Re\phi_{32}&-\Re\phi_{32}
&0&0\\

0&0 &4(k-1)\phi_{22}&
\begin{matrix}
-2k^{2}+3k-2-2\phi_{22}^{2}\\+3\phi_{33}
+4\Im\phi_{42}
\end{matrix}
&0
&-3\Im\phi_{32}&-\Im\phi_{32}
&0&0\\

0
&
0
&
\begin{matrix}
 8\Re\phi_{43} -8\phi_{22}\Re\phi_{32} \\
 -14(k-1)\Im \phi_{32}
\end{matrix}
&
\begin{matrix}
 8\Im\phi_{43}-8\phi_{22}\Im\phi_{32}\\
 +2(k-1)\Re \phi_{32}
\end{matrix}
&
0
&
\begin{matrix}
\frac{k(2k+3)}3\\-\phi_{33}
\end{matrix}
&
\begin{matrix}
6(k-1)\\
\times\phi_{22}
\end{matrix}
&
0
&0
\end{pmatrix}.
$$
\normalsize

Calculating the determinant we obtain
$\det A=\frac14(k-1) \det B$,
where

\begin{equation}\Label{B}
B=
\begin{pmatrix}
-6\Re\phi_{32}&6\Im\phi_{32}	&-2\phi_{22}&2(k-1)\\

\begin{matrix}
2k^{2}-3k+2+2\phi_{22}^{2}\\
-3\phi_{33}
-4\Re\phi_{42}
\end{matrix}
&
4(k-1)\phi_{22}
&-3\Re\phi_{32}&-\Re\phi_{32}
\\

4(k-1)\phi_{22}&
\begin{matrix}
-2k^{2}+3k-2-2\phi_{22}^{2}\\+3\phi_{33}
+4\Im\phi_{42}
\end{matrix}
&-3\Im\phi_{32}&-\Im\phi_{32}
\\

\begin{matrix}
 8\Re\phi_{43} -8\phi_{22}\Re\phi_{32} \\
 -14(k-1)\Im \phi_{32}
\end{matrix}
&
\begin{matrix}
 8\Im\phi_{43}-8\phi_{22}\Im\phi_{32}\\
 +2(k-1)\Re \phi_{32}
\end{matrix}
&
\begin{matrix}
\frac{k(2k+3)}3\\-\phi_{33}
\end{matrix}
&
\begin{matrix}
6(k-1)\\
\times\phi_{22}
\end{matrix}

\end{pmatrix}.
\end{equation}
The statement of the lemma now readily follows.
\end{proof}

\begin{Rem}
We note here that it is not necessary to require the full prenormalization
given in Proposition~\ref{transcoeffs} in order to guarantee that the expression
\eqref{B} for the matrix $B$ above gives rise to the characteristic polynomial. Indeed, it is enough to require that $\varphi'$ just satisfies \eqref{addinnorm'}, since in this case, \eqref{fl0} implies that \eqref{fl0'} holds; i.e., we must have $f(z,0) = z$ for any transformation respecting
the prenormalization \eqref{addinnorm'}, and hence $\phi'_{ab}=\phi_{ab}$.
\end{Rem}


If $M$ is in nonresonant, as described in Definition \ref{generalpos}, then it follows from \eqref{reimmatrix} that we can inductively require the following additional normalization for $k\geq 2$:
\begin{equation}\Label{completenorm}
\phi_{10k}=\phi_{11k}=\phi_{21k}=\phi_{22k}=\phi_{32k}=\phi_{33k}=0,
\end{equation}
which will completely determine the variables \eqref{transvars}, i.e.,
$$
f_{0,k-1},\ f_{1,k-1},\ f_{2,k-1},\ \Re g_{0k},\ g_{1k},
$$
in Proposition \ref{transcoeffs}. It follows from \eqref{reimmatrix}, and a straighforward induction on $k\geq 2$ using also Proposition \ref{transcoeffs}, that the only transformation preserving the complete normalization described above is the identity mapping.
More generally, if $M$ has resonances $k$, we can still obtain the equations
\eqref{completenorm} for all non-resonant $k$.
We summarize this result in the following theorem.

\begin{Thm}\Label{main1} Let $M$ be a formal hypersurface through $0$ in $\bC^2$, satisfying the assumptions described in Subsection $\ref{setup}$. Assume furthermore that $M$ is in general position at $0$. Then there are formal holomorphic coordinates $(z,w)$ at $0$ such that $M$ is given as a formal graph
\begin{equation}
\Im w=\phi(z,\bar z,\Re g),
\end{equation}
where the formal (Hermitian) power series $\phi(z,\bar z,u)$ is of the form
\begin{equation}
\phi(z,\bar z,u)=\sum_{a,b\geq 0}\phi_{ab}z^a\bar z^bu+\sum_{a,b\geq 0\atop k\geq 2}\phi_{abc}z^a\bar z^bu^k
\end{equation}
satisfying the following normalization conditions
\begin{equation}\Label{normalform1}
\phi_{11}=1,\ \phi_{a0}=\phi_{l1}=\phi_{a0k}=\phi_{l+1,1k}=0,\quad a\geq 0,\, k, l\geq 2\,
\end{equation}
and
\begin{equation}\Label{normalform2}
\phi_{10k}=\phi_{11k}=\phi_{21k}=\phi_{22k}=\phi_{32k}=\phi_{33k}=0,\quad k\geq 2.
\end{equation}
Moreover, the only formal transformation of the form
\begin{equation}\Label{transagain}
z'=z+f(z,w),\quad w'=w+g(z,w),
\end{equation}
where $f$ and $g$ are formal holomorphic power series with $f(0,0)=g(0,0)=f_z(0,0)=g_w(0,0)=0$, that preserves the normalization \eqref{normalform1} and \eqref{normalform2} is the identity, i.e., $f\equiv g\equiv 0$.

Furthermore, without assuming $M$ to be in general position, we still obtain its formal normalization given by all equations
\eqref{normalform1} and those in \eqref{normalform2} for all non-resonant $k$.
\end{Thm}

\begin{Rem}{\rm We note that there is some redundancy in the conditions \eqref{normalform1} and \eqref{normalform2}. The reason we present the conditions in this way here is so that the reader can keep track of which conditions come from the prenormalizations in Subsection \ref{prenormsec} (those in \eqref{normalform1}) and which come from the final normalization in Subsection \ref{completenormsec} (those in \eqref{normalform2}). In Theorem \ref{main0}, we have eliminated this duplication of conditions, and present the results in a form that closely mimics the Chern-Moser normal form.
}
\end{Rem}

To round out the discussion, we note that a general invertible transformation
$$(z',w')=(F(z'',w''),G(z'',w''))$$
 preserving the normalization in Theorem \ref{main1} can be factored as $(z,w)=(\alpha z'',s w'')$ composed with a transformation of the form \eqref{transagain}; in order to preserve the real tangent space to $M$ at $0$, we need to require $s\in \bR^*:=\bR\setminus \{0\}$, and in order to preserve $\phi_{11}=1$, we must require $|\alpha|=1$. Since the linear transformation $(z,w)=(\alpha z'',s w'')$ preserves the normalization, we conclude that the group $G:=S^1\times \bR^*$ acts on the space of normal forms and the isotropy group of $M$ at $0$ is a subgroup of $G$. Moreover, the uniqueness part of Theorem \ref{main1} implies the following: {\it Any formal holomorphic transformation that preserves the normal form in Theorem $\ref{main1}$ is of the form $(z,w)\mapsto (\alpha z,sw)$ with $(\alpha,s)\in S^1\times \bR^*$.} Theorem \ref{main0} now follows easily by writing the defining equation of $M$ in the form
$$
\Im w= \Re w\left(|z|^2+\sum_{a,b\geq 0} N_{ab}(\Re w)z^a\bar z^b\right),
$$
and translating the conditions in Theorem \ref{main1} into conditions on $N_{ab}(u)$.

\section{Examples}

We conclude this paper by giving a few examples.

\begin{Exa} Consider a hypersurface $M\subset\bC^2$ of the following form
\begin{equation}
\Im w=  \phi(z,\bar z,\Re w), \quad \phi(z,\bar z,u) = u |z|^2+ u^2\psi(z,\bar z,u),
\end{equation}
where $\psi(z,\bar z,u)$ is such that $\phi(z,\bar z,u)$ satisfies
\eqref{orgnorm2'} and \eqref{addinnorm'},
for example,
$$\psi(z,\bar z,u)=\theta(|z|^2,u)$$
where $\theta(x,u)$ satisfies $\theta_x(0,u)=0$.
In view of \eqref{fl0'} all terms involving $\phi_{ab}$ in
$\det B$ in \eqref{B} are the same as in the normal form, and hence are 0, and we compute
$$
\det B=\frac{2}{3}k(2k+3)(k-1)(2k^2-3k+2)^2.
$$
Since the roots of $2k^2-3k+2$ are not real, we conclude that $M$ is nonresonant at $0$. Therefore, we can put $M$ into normal form as described in Theorem \ref{main0}, i.e., eliminate terms of the form $|z|^4 u^k$, $z^3\bar z^2u^k$, and $|z|^6u^k$ in $\psi(z,\bar z,u)$). The stability group of $M$ is a subgroup of $S^1$, unless $\psi$ after normalization vanishes completely.
\end{Exa}

\begin{Exa}\Label{GenEx}
If $M$ is given by an equation of the form
\[ \Im w = \Re w \left(|z|^2 + \frac{C}{4} |z|^4 + \frac{D}{36} |z|^6 + O(|z|^8)\right) +
O((\Re w)^2),\]
then the characteristic polynomial is given by
\begin{equation}\Label{GE}
-8 i (k-1) \left|24 (k-1)^2+6 i C (k-1)+3 C^2-D+12\right|^2 \left(48 (k-1)^2+27 C^2-8 D+96\right).
\end{equation}
The first two factors do not have any integral roots $k\ge 2$ provided that
$C\neq 0$. In this case, there is for any integer $k\geq 2$ an unique
$D$ such that the characteristic polynomial has exactly
that resonance $k$. If on the other hand, $C=0$, the characteristic
polynomial is given by
\[64 i (k-1) \left(D-24 k^2+48 k-36\right)^2 \left(D-6 \left(k^2-2
   k+3\right)\right).\]
The reader can easily check that the last two factors have no real roots if $D< -12$, one root each at $k=1$ if $D=-12$, and two (distinct) roots each, symmetric about about $k=1$ when $D>-12$. In particular, for $C=0$, the hypersurface $M$ has either zero or two resonances.
\end{Exa}

We note that if $M$ satisfies Conditions (1') and (2') in the introduction, but has a resonance at $k=k_0\geq 2$, then we cannot in general achieve the normalization \eqref{normalform2} at $k=k_0$. We can make a choice of the derivatives of $f$ and $g$ in \eqref{transvars} at $k=k_0$, and then proceed with the normalization for $k>k_0$ (until the next resonance, if it exists). However, the choice of \eqref{transvars} at $k=k_0$ will in general affect the corresponding coefficients \eqref{normalform2}. Therefore, the existence of a resonance $k=k_0$ does not necessarily imply that the derivatives \eqref{transvars} at $k=k_0$ of an automorphism of $M$ is not determined by previous ones. There are, however, known examples of $M\subset \bC^2$ satisfying Conditions (1') and (2'), whose stability groups at $0$ are not determined by 1-jets (see \cite{Kow03}, \cite{Zaitsev02}, \cite{KL14}). Such hypersurfaces cannot be in general position (i.e., must have resonances) at 0 by Corollary \ref{maincor}, and the failure of $1$-jet determination is caused by the resonances. We mention two such examples here (of the form in Example \ref{GenEx}), where a resonance $k=k_0$ actually corresponds to indeterminacy of the derivatives \eqref{transvars} for $k=k_0$ in automorphisms of $M$; more examples can be found in the list in \cite{KL14}.

\begin{Exa} For a positive integer $m$, consider the following $M_m\subset \bC^2$,
\begin{equation}
\Im w=i\Re w \frac{1-q_m(2m |z|^2)}{1+q_m(2m|z|^2)} = \Re w \left( |z|^2 + \left(\frac{2 m^2}{3}+\frac{1}{3}\right) |z|^6 + \dots  \right) ,
\end{equation}
where
\begin{equation}
q_m(x)=e^{(i/m)\arcsin x}.
\end{equation}
It is readily checked that $M_m$ satisfies (1') and (2') at $0$, and comparing with the formula for the characteristic polynomial in \eqref{GE} (with $C=0$), we see that it is given by
\[ -221184 i (k-1) \left((k-1)^2-4 m^2\right) \left((k-1)^2-m^2\right)^2. \]
Its resonances are therefore given by $k = m+1$ and $k=2m+1$.
We note (cf.\ \cite{Zaitsev02}, \cite{KL14}) that the following local 1-parameter family of biholomorphisms belong to its stability group at $0$:
\begin{equation}
H_t(z,w)=\left(\frac{z}{(1-tw^{2m})^{1/2}},\frac{w}{(1-tw^{2m})^{1/2m}}\right),\quad t\in \bR.
\end{equation}
We note that the jets $j^{2m}_0 H_t$ agree for all $t\in \bR$, but the derivatives \eqref{transvars} for $k=2m+1$ depend on the parameter $t$.
\end{Exa}

\begin{Exa} The following example (corresponding to $C \neq 0$ in Example \ref{GenEx})
illustrates that even a single resonance can be responsible for the presence of automorphisms not determined by their $1$-jets.
We let  $q_T (|z|^2)$ be the unique solution to
\[ u q_T'(u) = \frac{\tan(q_T (u))}{1 + T \tan(q_T(u))}, \quad q_T(0)=0,\ q'_T(0)=1, \]
where $T\in \R\setminus\{0\}$.
Then, for any $m\in\N$, the hypersurface $M_{m,T}$ defined by
\begin{equation*}\Im w = \Re w \tan \left( \frac{q_T (m |z|^2)}{m}\right) = \Re w
 \left( |z|^2 -  m T |z|^4 +  \left(\frac{2+ m^2(9 T^2+1)}{6}\right)|z|^6 +\dots  \right)
\end{equation*}
also satisfies (1') and (2') at $0$. It
has an infinitesimal CR automorphism given by
\begin{equation}\Label{Mauto}
X = \frac{1}{m} \left(\frac{1}{2} + i T\right) zw^m \frac{\partial}{\partial z} + w^{m+1}\frac{\partial}{\partial w},
\end{equation}
as can be seen from a computation carried out in \cite[Lemma 10 and 4]{KL14}.
The resonances of its characteristic polynomial, by the observation that the first two factors of the characteristic polynomial in \eqref{GE} do not have any if $T\neq 0$, are the integral roots $k\geq 2$ of the polynomial
\[ k^2 - 2k + 1 - m^2; \]
that is, only $k=m+1$ occurs. The infinitesimal CR automorphism \eqref{Mauto}  illustrates the failure of the conclusion in
Corollary~\ref{maincor} in this case.
\end{Exa}


\def\cprime{$'$}


\begin{thebibliography}{BER99}

\bibitem[BER00]{BER00}
M.~S. Baouendi, Peter Ebenfelt, and Linda~Preiss Rothschild.
\newblock Convergence and finite determination of formal {CR} mappings.
\newblock {\em J. Amer. Math. Soc.}, 13:697--723, 2000.

\bibitem[BER99]{BER99a}
M.~Salah Baouendi, Peter Ebenfelt, and Linda~Preiss Rothschild.
\newblock {\em Real submanifolds in complex space and their mappings},
  volume~47 of {\em Princeton Mathematical Series}.
\newblock Princeton University Press, Princeton, NJ, 1999.

\bibitem[BMR02]{BMR02}
M.~Salah Baouendi, Nordine Mir, and Linda~Preiss Rothschild.
\newblock Reflection ideals and mappings between generic submanifolds in
  complex space.
\newblock {\em Journal of Geometric Analysis}, 12:543--580, 2002.
\newblock 10.1007/BF02930653.

\bibitem[CM74]{CM74}
S.~S. Chern and J.~K. Moser.
\newblock Real hypersurfaces in complex manifolds.
\newblock {\em Acta Math.}, 133:219--271, 1974.

\bibitem[ELZ03]{ELZ03}
P.~Ebenfelt, B.~Lamel, and D.~Zaitsev.
\newblock Finite jet determination of local analytic {CR} automorphisms and
  their parametrization by 2-jets in the finite type case.
\newblock {\em Geom. Funct. Anal.}, 13(3):546--573, 2003.

\bibitem[E98a]{E98a}
Peter Ebenfelt.
\newblock New invariant tensors in {CR} structures and a normal form for real
  hypersurfaces at a generic {L}evi degeneracy.
\newblock {\em J. Differential Geom.}, 50(2):207--247, 1998.

\bibitem[E98b]{E98b}
Peter Ebenfelt.
\newblock Normal forms and biholomorphic equivalence of real hypersurfaces in
  {$\bold C^3$}.
\newblock {\em Indiana Univ. Math. J.}, 47(2):311--366, 1998.

\bibitem[EHZ05]{EHZ05}
Peter Ebenfelt, Xiaojun Huang, and Dmitri Zaitsev.
\newblock The equivalence problem and rigidity for hypersurfaces embedded into
  hyperquadrics.
\newblock {\em Amer. J. Math.}, 127(1):169--191, 2005.

\bibitem[ELZ09]{ELZ09}
Peter Ebenfelt, Bernhard Lamel, and Dmitri Zaitsev.
\newblock Degenerate real hypersurfaces in $\mathbb c^2$ with few
  automorphisms.
\newblock {\em Trans. Am. Math. Soc.}, 361(6):3241--3267, 2009.

\bibitem[J07]{Juhlin07}
Robert Juhlin.
\newblock Convergence of formal cr mappings at points of 1-infinite type.
\newblock {\em preprint}, 2007.

\bibitem[JL13]{JuhlinLamel13}
Robert Juhlin and Bernhard Lamel.
\newblock On maps between nonminimal hypersurfaces.
\newblock {\em Math. Z.}, 273(1-2):515--537, 2013.

\bibitem[Kol05]{Kolar05}
Martin Kol{\'a}{\v{r}}.
\newblock Normal forms for hypersurfaces of finite type in {${\mathbb{C}}^2$}.
\newblock {\em Math. Res. Lett.}, 12(5-6):897--910, 2005.

\bibitem[Kol12]{Kolar12}
Martin Kol{\'a}{\v{r}}.
\newblock Finite type hypersurfaces with divergent normal form.
\newblock {\em Math. Ann.}, 354(3):813--825, 2012.

\bibitem[KL14]{KL14}
Martin Kol{\'a}{\v{r}} and Bernhard Lamel.
\newblock Holomorphic equivalence and nonlinear symmetries of ruled
   hypersurfaces in $\Bbb{C}^2$
\newblock {\em J. Geom. Anal.}, 25(2):1240--1281, 2015.

\bibitem[KMZ14]{KMZ14}
Martin Kolar, Francine Meylan, and Dmitri Zaitsev.
\newblock Chern-{M}oser operators and polynomial models in {CR} geometry.
\newblock {\em Adv. Math.}, 263:321--356, 2014.

\bibitem[KoL14]{KoL14}
Ilya Kossovskiy, Bernhard Lamel.
\newblock On the analyticity of CR-diffeomorphisms. arXiv.org.
\newblock{{\em {\tt http://arxiv.org/abs/1408.6711v1}}}, 2014.

\bibitem[KoL15]{KoL15}
Ilya Kossovskiy, Bernhard Lamel.
\newblock New extension phenomena for solutions of tangential Cauchy-Riemann Equations. arXiv.org.
\newblock {{\em \tt http://arxiv.org/abs/1507.06151v1}}, 2015.

\bibitem[KZ14a]{KoZ14a}
Ilya Kossovskiy, Dmitri Zaitsev.
\newblock Convergent normal form for real hypersurfaces at generic levi
  degeneracy.
\newblock {\em {\tt http://front.math.ucdavis.edu/1405.1743}}, 2014.

\bibitem[KZ14b]{KoZ14b}
Ilya Kossovskiy, Dmitri Zaitsev.
\newblock Convergent normal form and canonical connection for hypersurfaces of finite type in $\mathbb C^2$.
\newblock {\em {\tt http://front.math.ucdavis.edu/1409.7506}}, 2014.

\bibitem[Kow02]{Kow03}
R.~Travis Kowalski.
\newblock A hypersurface in {$\Bbb C^2$} whose stability group is not
  determined by 2-jets.
\newblock {\em Proc. Amer. Math. Soc.}, 130(12):3679--3686 (electronic), 2002.

\bibitem[Kow05]{Kow05}
R.~Travis Kowalski.
\newblock Rational jet dependence of formal equivalences between real-analytic
  hypersurfaces in {$\Bbb C^2$}.
\newblock {\em Pacific J. Math.}, 220(1):107--139, 2005.

\bibitem[S91]{Stanton91}
N.~Stanton.
\newblock A normal form for rigid hypersurfaces in $\bC^2$.
\newblock {\em Amer. J. Math.}, 113:877--910, 1991.

\bibitem[W82]{Wong82}
Philip~P. Wong.
\newblock A construction of normal forms for weakly pseudoconvex {CR} manifolds
  in {${\bf C}^{2}$}.
\newblock {\em Invent. Math.}, 69(2):311--329, 1982.

\bibitem[Z02]{Zaitsev02}
Dmitri Zaitsev.
\newblock Unique determination of local {CR}-maps by their jets: a survey.
\newblock {\em Atti Accad. Naz. Lincei Cl. Sci. Fis. Mat. Natur. Rend. Lincei
  (9) Mat. Appl.}, 13(3-4):295--305, 2002.
\newblock Harmonic analysis on complex homogeneous domains and Lie groups
  (Rome, 2001).

\end{thebibliography}

\end{document}